# On the Preliminary Design of Multiple Gravity-Assist Trajectories

M. Vasile.[*]

*Politecnico of Milano, Milan, 20156, Italy*

P. De Pascale[†]

*University of Padova, Padova, 35131, Italy*

**In this paper the preliminary design of multiple gravity-assist trajectories is formulated as a global optimization problem. An analysis of the structure of the solution space reveals a strong multimodality, which is strictly dependent on the complexity of the model. On the other hand it is shown how an oversimplification could prevent finding potentially interesting solutions. A trajectory model, which represents a compromise between model completeness and optimization problem complexity is then presented. The exploration of the resulting solution space is performed through a novel global search approach, which hybridizes an evolutionary based algorithm with a systematic branching strategy. This approach allows an efficient exploration of complex solution domains by automatically balancing local convergence and global search. A number of difficult multiple gravity-assist trajectory design cases demonstrates the effectiveness of the proposed methodology.**

## Nomenclature

$D_l$     =     *lth* subdomain

$F$     =     generic objective function

$f$     =     fitness value

---

[*] Assistant Professor, Department of Aerospace Engineering, Via La Masa 34, 20156, Milano, Italy, AIAA Member.
[†] PhD candidate, Center for Space Studies CISAS, Via Venezia 15, 35131, Padova, Italy, AIAA Member.

| | | |
|---|---|---|
| $h_k$ | = | pericenter altitude, km |
| $k$ | = | phase identification index |
| $L_b$ | = | lower bound on the design parameters |
| $N_p$ | = | number of planets of a sequence |
| $p_k$ | = | planet identification number |
| $R$ | = | residual on the constraint |
| $S$ | = | individual exploration region |
| $T$ | = | time of flight, day |
| $t_0$ | = | launch date, MJD2000 |
| $t_{ds}$ | = | time of deep space maneuver, days |
| $U_b$ | = | upper bound on the design parameters |
| $v^\infty$ | = | hyperbolic excess velocity, km/s |
| $\mathbf{n}_\Pi$ | = | vector normal to the plane of the hyperbola |
| $\mathbf{q}$ | = | quaternion |
| $\mathbf{r}$ | = | heliocentric position in Cartesian coordinate, km |
| $\mathbf{u}_i$ | = | Euler vector for the rotation of the incoming velocity |
| $\tilde{\mathbf{v}}_i$ | = | flyby incoming relative velocity, km/s |
| $\tilde{\mathbf{v}}_o$ | = | flyby outgoing relative velocity, km/s |
| $\mathbf{y}$ | = | solution vector |
| $\alpha$ | = | right ascension of the launch asymptote, rad |
| $\delta$ | = | declination of the launch asymptote, rad |
| $\varepsilon$ | = | epoch of a deep space maneuver as a percentage of transfer time |
| $\gamma$ | = | deflection angle, rad |
| $\Delta v$ | = | change of velocity, km/s |
| $\eta_i$ | = | plane rotation angle, rad |
| $\rho$ | = | migration radius |
| $\sigma$ | = | weighting factor |

$\varphi_{Dq}$  =  subdomain fitness

$\psi_{Dq}$  =  subdomain qualification index

$\omega_{Dq}$  =  subdomain density

## I. Introduction

THE access to high Δv targets in the Solar System has been made possible by the enabling concept of gravity assist,[1,2] through which many missions with challenging objectives were successfully accomplished with significant reductions in launch energy and in transfer time. Past missions, like Mariner[3] 10, Voyager[4] 1 and 2, Pioneer 10 and 11 and more recent missions like Galileo,[5] Ulysses and Cassini[6] made extensively use of gravity-assist maneuvers; future missions to completely or partially unexplored planets, such as Mercury[7,8] or Pluto,[9,10] will be possible thanks to gravity-assist maneuvers.

The availability of new propulsion systems, such as low-thrust engines, while widening the possibility of flights in the Solar System, has not diminished the importance of gravity-assist maneuvers. This is true at least for two different reasons: many missions will still resort to pure ballistic trajectories using only impulsive corrections; optimal gravity-assist sequences are important in the design of low-thrust gravity-assist trajectories[7].

Multiple gravity-assist trajectories have been extensively investigated over the last forty years, and their preliminary design has been approached mainly relying on the experience of mission analysts and following simplifying assumptions in unison with systematic searches or some simple analysis tools such as the Tisserand's graph.[1,10] This graphical tool, however, turns out to be unsuitable for the investigation of transfers to asteroids and comets with highly eccentric or inclined orbits, and for the analysis of transfer options that require additional corrective maneuvers, performed by either high-thrust or low-thrust engines.

On the other hand, since at the early stage of the design of a space mission a number of different options is generally required, it would be desirable to automatically generate many optimal or nearly optimal solutions over the range of the design parameters (escape velocity, launch date, time of flight, etc…), accurately enough to allow a correct trade-off analysis. This could translate into the need for an exhaustive characterization of highly complex solution domains deriving from sophisticated trajectory models.

Recently, different attempts have been carried out toward the definition of automatic design tools, although so far most of these tools have been based on systematic search engines. An example is represented by the automatic tool for the investigation of multiple gravity-assist transfers, called STOUR, originally developed by JPL[5] and subsequently enhanced by Longuski et al. at Purdue university.[11,12,13,14] This tool has been extensively used for the preliminary investigation of interplanetary trajectories to Jupiter and Pluto,[11,12] for the design of the tour of Jovian moons[15] and for Earth-Mars cycling trajectories.[16]

It should be noted that such a systematic search methodology although theoretically complete, could require long computational times, in order to identify all the suitable options. Although completeness is a desirable property of any strategy that aims at finding all optimal or suboptimal options, it is also true that this property could remain merely theoretical if the complexity of the problem makes it practically intractable by means of a systematic approach (NP-hard problems).

In the last ten years, different forms of stochastic search methods have also been applied to orbit design, starting from the work of Coverstone[17] et al. on the use of multi-objective genetic algorithms for the generation of first guess solutions for low-thrust trajectories, to more recent works on the use of single-objective genetic algorithms for ballistic transfers[18] or to the use of hybrid evolutionary search method for preliminary design of weak stability boundaries (WSB) and interplanetary transfers.[19,20]

In this paper we investigate the possibility of using stochastic based methods for the solution of multiple gravity-assist trajectories and we propose a novel approach for the preliminary automatic design of complex trajectories. At first, an analysis of the relation between model and problem complexity will be presented. This will highlight how the structure of the solution space is strongly dependent on the mathematical modeling of gravity-assist maneuvers and transfer trajectories. In particular, it will be shown that problem complexity grows significantly if deep space maneuvers are included and physical constraints on gravity maneuvers are satisfied. On the other hand, models with no deep space maneuvers but with $\Delta v$ matching at the swing-by planet have a simpler structure[21] and the related optimization problem can be efficiently solved with a clever systematic search[22] (enumerative, multi-grid or branch and prune). This is even truer if the $\Delta v$ matching can be performed with an aero-gravity maneuver instead of a propelled maneuver.[23]

In order to handle the complex structure of more realistic models, in this paper the use of a novel global search algorithm based on a hybridization of an evolutionary based search and a deterministic branching strategy[19,20] is proposed. An important feature of this approach is that it combines a stochastic based search with a complete deterministic decomposition of the solution space. This allows a loose specification of the bounds on the initial solution domain and a parallel exploration of different promising areas of the search space. The output is therefore an exhaustive characterization of the search space and a number of alternative families of solutions.

This hybrid approach is here applied to the solution of a number of complex interplanetary transfers both to planets and to comets and asteroids, showing its effectiveness and providing some new solutions to known mission analysis problems.

## II. Trajectory Model Analysis

An engineering design problem can always be tackled with a two-stage approach: problem modeling and problem solution, where often the search for a solution is represented by an optimization procedure. Modeling is the task of transcribing a physical phenomenon into a mathematical representation. The modeling stage has a particular influence on the definition and development of preliminary design methodologies since there is always a trade-off between the precision of the required solution and the computational cost associated to its search:[20] different models intrinsically contain different kinds of solutions and can favor or not their identification.

This issue becomes of significant importance when a large number of good first guess solutions has to be efficiently generated for an exhaustive preliminary assessment of complex engineering problems. In this case efficiency is quantified as the ratio between the number of useful solutions and associated computational time. These considerations obviously also apply to trajectory design. Therefore the two above-mentioned stages, which are in fact mutually dependent, must be properly defined during the development of an effective design tool for the preliminary investigation of complex interplanetary transfers.

In particular, the modeling process requires the identification of the most important features of the trajectory that will be analyzed and must reproduce the completeness of the problem under investigation, while reducing its complexity. This is a trivial consideration, that has non trivial consequences on the effectiveness of the design phase, since an oversimplified model could lead to the loss of interesting solutions.

On the other hand a proper mathematical model, which accurately reproduces a physical phenomenon, is likely to require more efficient search methods, in order to find a specific solution. Therefore the problem solution stage needs proper search mechanisms or approaches that allow the identification of all relevant solutions in a given solution domain. This raises to the additional issue of the completeness of the search: if the problem is at least NP-hard, a complete search may not be practically possible since the number of function evaluations to prove the optimality of a solution could grow exponentially with problem dimension.

In the following, the attention will be focused on some simple trajectory models of increasing complexity in order to derive a good compromise between computational cost and solution accuracy. Considering the typical multiple gravity-assist trajectories that have been designed and flown so far, some general features can be considered relevant in order to maintain the required richness of the search space and to identify all the families of solutions that could be potentially interesting for the design of an interplanetary mission. In particular a full 3D model both for the trajectory and for the gravity-assist maneuvers has been developed including deep space maneuvers (DSM) and using the analytical ephemeris of celestial bodies. The benefits of such a modeling approach can be seen in the design of missions to Pluto, Mercury or to the Sun, which require consideration of the real inclination of the orbit of the planet or of the final heliocentric orbit, and in the design of missions to near earth objects. This particular choice is compared, in terms of search space complexity, to a simpler model in which DSM are neglected. This simple modification prevents consideration of some classes of interesting solutions such as, resonant or almost resonant swing-by or free orbits before the encounter with a celestial body.

Since the physics of the Solar System allows us the adoption of a patched-conic approximation of a multiple gravity-assist trajectory, a complete transfer trajectory can be reduced to the sum of a number of smaller sub-problems with a finite number of design variables. As will be shown in the following, each sub-problem may produce complex search domains, typically non-convex and multimodal.

**A. Two- and Three-Impulse Transfers**

A simple two-impulse direct transfer from a planet $P_1$ to a planet $P_2$ can be modeled as a function of the departure date $t_0$ and of the arrival date, expressed as the sum of the departure date and of the time of flight $T$. The two required impulsive maneuvers at departure $\Delta\mathbf{v}_1$ and at arrival $\Delta\mathbf{v}_2$ can be computed solving a Lambert's problem[2] from $P_1$ to $P_2$ in a given time $T$.

An extension of this simple two-impulse model can be obtained by propagating analytically the initial state, given by the position of $P_1$ and by its velocity plus $\Delta\mathbf{v}_1$, for a time $t_{ds}=\varepsilon T$, up to a point $M_1$ and then by solving a Lambert's problem from $M_1$ to $P_2$ as in Fig. 1. The resulting discontinuity $\Delta\mathbf{v}_s$ in the velocity at $M_1$ represents a deep space maneuver or a correction on the direct transfer. The total cost of the transfer can be expressed as:

$$f(t_0, \varepsilon, T) = \Delta v_1 + \Delta v_s + \Delta v_2 \quad (1)$$

The structure of the solution space for a two-impulse transfer is directly related to the synodic period of the two bodies $P_1$ and $P_2$ and to the orbital elements of their respective orbits, as can be seen in Fig.2a where all solutions for a direct Earth-Mars transfer with a total $\Delta v$ lower than 15 km/s has been represented. On the other hand, a deep space maneuver, increases the number of feasible paths from $P_1$ to $P_2$ as can be seen in Fig. 2b where the function in Eq. (1) has been plotted for all the optimal values of $\varepsilon$, keeping the $\Delta v$ at departure aligned along the velocity of the Earth.

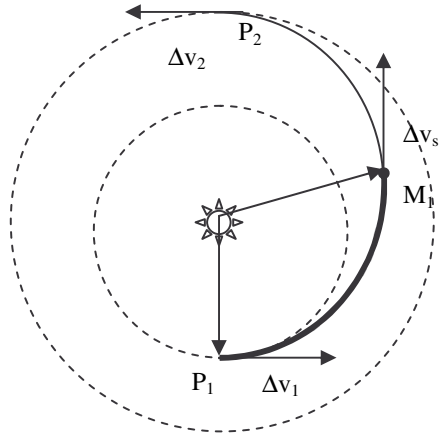

**Fig. 1 Schematic of a three-impulse transfer.**

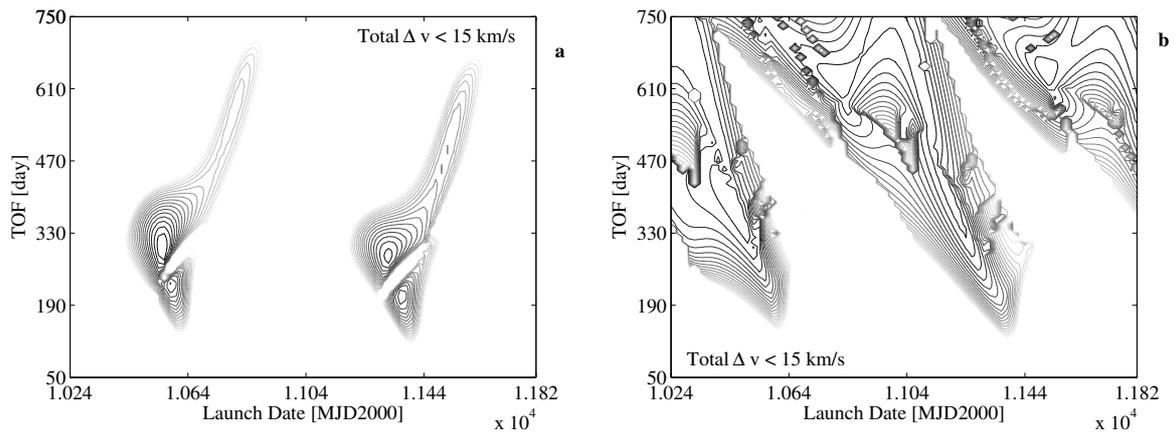

**Fig. 2 Solution space for: a) a two-impulse transfer, b) a three-impulse transfer.**

**B. Gravity-Assist Maneuvers**

Gravity-assist maneuvers are here modeled in three dimensions using a linked conic approximation: the maneuver is instantaneous and the sphere of influence is collapsed into a single point. Other two conditions define the physical properties of the gravity maneuver: the incoming and outgoing relative velocities $\tilde{v}_i$ and $\tilde{v}_o$ must have equal modulus and the deflection angle is a function of the pericenter altitude $h$ and of the incoming relative velocity. These conditions typically translate into the following set of constraint equations on position and velocity:

$$\mathbf{r}_i = \mathbf{r}_o = \mathbf{r}_p, \tag{2}$$

$$\tilde{v}_i = \tilde{v}_o, \tag{3}$$

$$\tilde{\mathbf{v}}_o^T \tilde{\mathbf{v}}_i = -\cos(2\beta(h))\tilde{v}_i^2 \tag{4}$$

The constraints on the incoming and outgoing position vectors $\mathbf{r}_i$ and $\mathbf{r}_o$ and on the planet position vector $\mathbf{r}_p$ in Eq. (2), can be explicitly solved, while the two conditions in Eq. (3),(4) on the relative velocities would require the solution of a constraint satisfaction problem. In order to avoid the introduction of these nonlinear constraints on GA maneuvers, a different model is used here. By the introduction of an auxiliary rotation angle $\eta$ it is possible to define the plane of the hyperbola where the rotation of the relative velocity occurs. This is necessary since in the linked conic approximation the point where the interplanetary trajectories pierces the sphere of influence is not defined and thus the plane where the planetocentric hyperbolic motion occurs remains undetermined. This auxiliary angle represents the rotation of a general plane $\Pi$, defined through its normal vector $\mathbf{n}_\Pi$, around the Euler vector $\mathbf{u}_i = \tilde{\mathbf{v}}_i/\tilde{v}_i$. Now the vector $\mathbf{u}_i$ and the rotation $\eta$ can be used to define the quaternions:

$$\mathbf{q}_\eta = [\mathbf{u}_i \sin\eta, \cos\eta]^T \tag{5}$$

Then, any rotation around $\tilde{\mathbf{v}}_i$ of a general vector $\mathbf{n}_i$ normal to $\tilde{\mathbf{v}}_i$ can be expressed as:

$$\mathbf{n}_\Pi = \mathbf{Q}_\eta(\mathbf{u}_i)\mathbf{n}_i \tag{6}$$

where $\mathbf{Q}_\eta$ is a rotation matrix defined by the quaternion $\mathbf{q}_\eta$. The normal vector $\mathbf{n}_\Pi$ completely defines the orbit plane of the swing-by hyperbola and can be used to define the rotation of the incoming relative vector onto the outgoing relative vector:

$$\tilde{\mathbf{v}}_o = \mathbf{Q}(\mathbf{n}_\Pi)\tilde{\mathbf{v}}_i \tag{7}$$

where $\mathbf{Q}(\mathbf{n}_\Pi)$ is the rotation matrix defined by the quaternion $\mathbf{q}_\gamma = [\mathbf{n}_\Pi \sin\gamma, \cos\gamma]^T$, with $\gamma=\pi-2\beta$.

There are different possible choices for the vector $\mathbf{n}_i$, depending on which reference plane is adopted. A possible choice is to take this vector as the normal direction to the incoming velocity vector in the *xy* plane as shown in Fig. 3a. An alternative form could be to define this vector as the normal direction to the plane containing the relative incoming velocity and the planet's velocity $\mathbf{V}_p$ as shown in Fig. 3b.

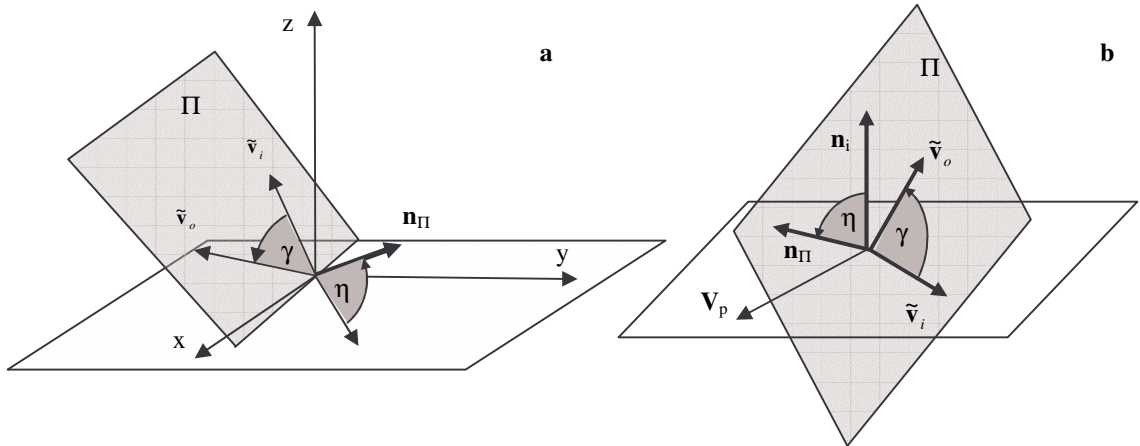

**Fig. 3 Schematic of the linked-conic model for a gravity assist; a) reference direction normal to the incoming velocity in the xy plane, b) reference direction normal to the relative and planet's velocity plane.**

Now, since a maneuver may be necessary in order to correct the post swing-by conditions, each gravity-assist maneuver can be associated to the following sub-problem: minimize the corrective deep space maneuver following the swing-by as a function of the swing-by characteristics and of the trajectory leg before the swing-by. The cost function of this sub-problem can be expressed as:

$$f(\mathbf{v}_i, t_i, h, \eta, \varepsilon, T) = \Delta v_s \tag{8}$$

If the incoming conditions $\mathbf{v}_i$, the swing-by epoch $t_i$ and the post swing-by transfer time $T$ are kept fixed, then $f$ can be plotted as a function of the angle $\eta$ and of the altitude of the swing-by $h$, for different values of parameter $\varepsilon$. Fig .4 show how the solution space progressively changes as the time of the deep space maneuver after the swing-by is delayed. As can be seen, due to the periodicity in $\eta$, the search domain is generally multimodal but an increase of $\varepsilon$ causes a significant increase in multimodality.

This behavior is similar for both the above-mentioned choices of the reference vector $\mathbf{n}_i$, although it has been noticed that choosing the normal direction to the incoming velocity in the *xy* plane leads to a more irregular solution space. Moreover this latter choice would make the discrimination of those flybys that maximize the energy gain, from those that maximize energy loss, quite difficult. In this respect, the former choice is much more intuitive since all flybys maximizing the energy gain have $\eta \in [\pi/2, 3/2\pi]$, while those which maximizing the energy loss have $\eta \in [-\pi/2, \pi/2]$. For this reason this choice has been adopted in the remainder of this work.

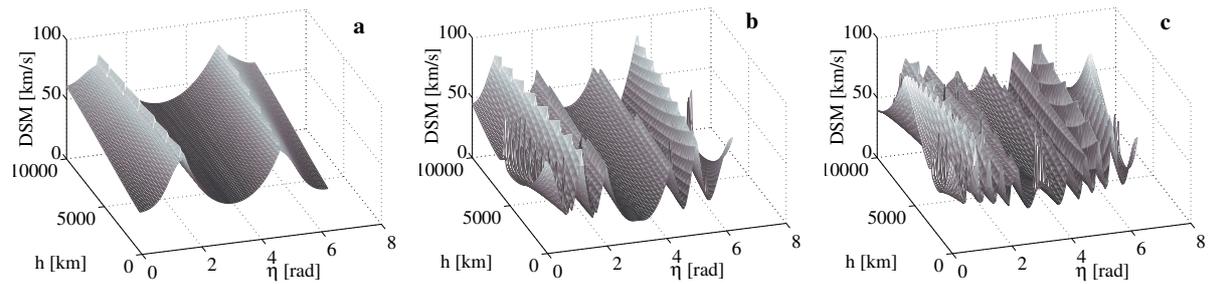

**Fig. 4 DSM cost after swingby for: a) impulsive maneuver at 10% of the transfer time, b) an impulsive maneuver at 37% the transfer time, c) an impulsive maneuver at 70% the transfer time.**

The consequences of the introduction of a deep space maneuver along a multiple gravity assist trajectory can be illustrated by the following example of a transfer from the Earth to Jupiter via a swingby of Venus. Two cases are analyzed: in the former case swingby model in Fig. 3b is used and a corrective deep space maneuver is performed after the swing-by while in the latter case a simple velocity matching at the swingby planet is performed (no deep space maneuver).

Position and velocity of the three planets are functions of the departure time $t_0$ and of the time of flight to Venus and to Jupiter $T_1$ and $T_2$. The impulsive maneuvers necessary to accomplish these transfers are $\Delta\mathbf{v}_1$ to leave the Earth, $\Delta\mathbf{v}_s$, which is either the required correction to match the outgoing velocity at the flyby in the case of no deep space maneuver, or the magnitude of the deep space maneuver itself, and $\Delta\mathbf{v}_2$ the arrival velocity at Jupiter. The total cost of the transfer can be expressed as:

$$f(t_0, T_1, h, \eta, \varepsilon, T_2) = \Delta v_1 + \Delta v_s + \Delta v_2 \tag{9}$$

with $\varepsilon$ being fixed to zero in the case of no deep space maneuver or free to vary between the instant of the flyby and 90% of the transfer time $T_2$ when a deep space maneuver is included.

The analysis of cost function in Eq. (9) has been performed by taking a set of 10000 randomly generated sample points and performing from each one a local search. The samples were generated using a Latin Hypercube distribution, and the Matlab function, fmincon was used for the local search setting the tolerance on optimality to 1e-6. This approach is considered satisfactory to give a rough characterization of the solution domain.

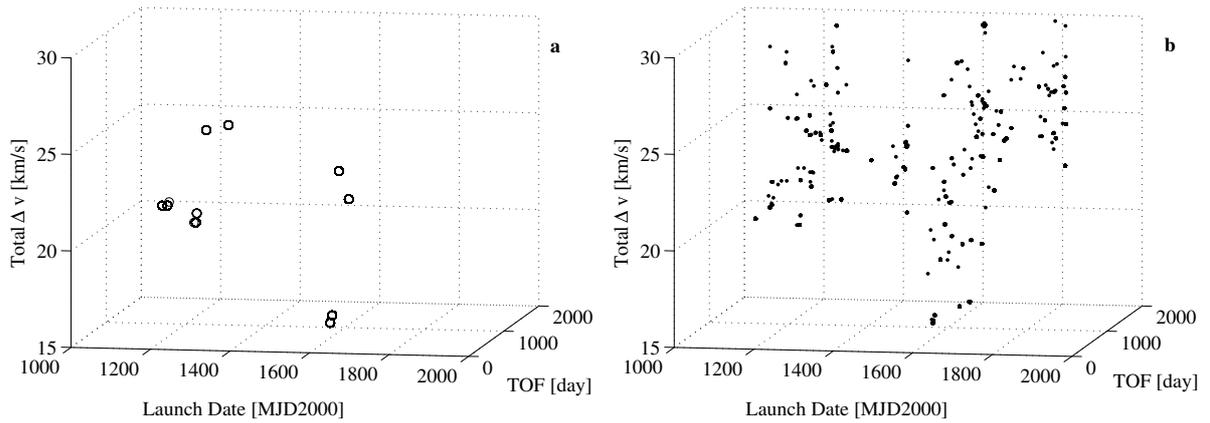

**Fig. 5 Distribution of the local minima; a) without deep space maneuver, b) with deep space maneuver.**

As can be seen from the comparison of Fig.5a and Fig.5b, where all the solutions with a total Δv lower than 30 km/s have been plotted, in the case of no deep space maneuver, there are two limited groups of local minima that are distributed over the launch dates corresponding to the synodic period of the Earth-Venus system. The situation becomes completely different in the case of a deep space maneuver performed after the swing-by of Venus. Although the globally optimal launch window is still the same, the number of opportunities for that window has increased and the distribution of the local minima is now almost continuous over the launch dates up to the epoch 1800 MJD. While an increase in the degrees of freedom due to the deep space maneuver could be somehow expected, less obvious is the benefit in terms of number of optimal solutions resulting from the increased complexity and multimodality of the problem.

The result of this simple test essentially shows that a simple model might not contain the required solutions but on the other hand a complex model could have such solutions so nested that finding them would be extremely difficult. In the following we will opt for a compromise between model complexity and model fidelity, relying on a global search algorithm for the analysis of the search space. This algorithm implements a general purposes strategy and is not endowed with additional information about the specific problem under study.

### III.     Complete Trajectory Model and Problem Formulation

Each single sub-problem, introduced and analyzed in the previous chapters, has been assembled into a complete model in which the trajectory is divided into a number of phases connecting a sequence of celestial bodies (the full trajectory model is reported in Fig.6). Given a sequence of $N_P$ planets, there exist $k=1,..., N_P-1$ phases, each of them beginning and ending with an encounter with a planet. Each phase $k$ is made of two conic arcs, the first ending where the second begins and having a discontinuity in the absolute heliocentric velocity at their matching point $M_k$. Given the transfer time $T_k$, relative to each phase $k$ and the variable $\varepsilon_k=[0,1]$, the matching point is then at $t_{dsk}=t_{k-1}+\varepsilon_k T_k$. The velocity vector at infinity, for the zero sphere-of-influence model, can be treated explicitly as a design parameter through the following definition:

$$\mathbf{v}_0^\infty = \sqrt{C_3}\left(\sin\delta\cos\alpha, \sin\delta\sin\alpha, \cos\delta\right) \tag{10}$$

with the angles $\delta$ and $\alpha$ respectively representing the declination and the right ascension of the escape asymptote. This choice allows to easily bound the escape velocity and asymptote direction within lower and upper values while adding the possibility of having a deep space maneuver in the first arc after the launch. This is often the case when escape velocity must be fixed due to the launcher capability or to the requirement of a resonant flyby of the Earth. Alternatively it is possible to use a simplified model in which the first leg from $P_0$ to $P_1$ is a simple Lambert's arc with no deep space maneuver. In this case the number of optimization parameters and degrees of freedom are reduced and the model is suitable for the assessment of sequences that fly directly to another planet after launch.

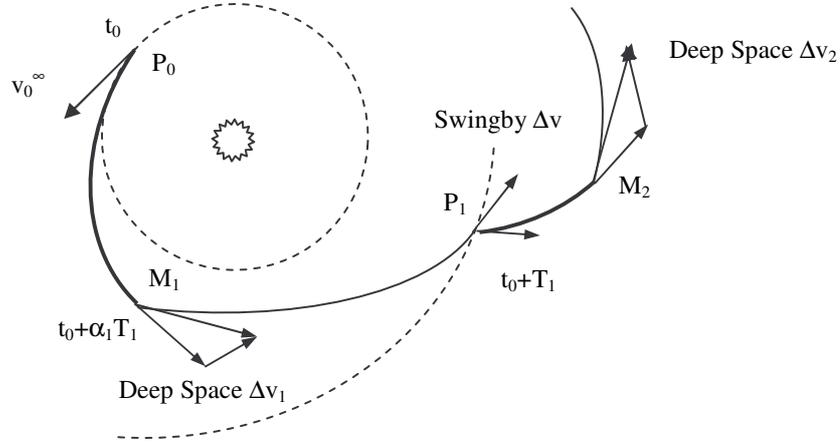

**Fig. 6 Schematic of multiple gravity-assist trajectories.**

Once the heliocentric initial velocity, which can be the result of a flyby maneuver or the asymptotic velocity after launch, is defined for each phase *k*, the trajectory is analytically propagated until time $t_{dsk}$. The second arc of phase *k* is then solved through a Lambert's algorithm, from $\mathbf{r}_{dsk}$, the Cartesian position of the deep space maneuver, to $\mathbf{r}_{k+1}$, the position of the target planet of phase *k*, for a time of flight $(1-\varepsilon_k)T_k$. Two subsequent phases are then joined together using the swing-by model given by Eq. (5) to Eq. (7). The solution vector for this model can be defined as:

$$\mathbf{y} = [v_0^\infty, \alpha, \delta, t_0, T_1, \varepsilon_1, \eta_1, h_1, ..., T_k, \varepsilon_k, \eta_k, h_k, ..., T_{Np-1}, \varepsilon_{Np-1}]^T \tag{11}$$

where $t_0$ is the departure date. Now, the design of a multigravity-assist transfer can be transcribed into a general nonlinear programming problem, with simple box constraints, of the form:

$$\begin{aligned} &\min F(\mathbf{y}) \\ &\text{with} \quad \mathbf{y} \in D \end{aligned} \tag{12}$$

One of the appealing aspects of this formulation is its solvability through a general global search method for box constrained problems, such as Evolutionary Algorithms (EA).[24] Depending on the kind of problem under study the objective function can be defined as:

$$F(\mathbf{y}) = \left\|\mathbf{v}_0^\infty\right\| + \sum_{k=1}^{N_p-1}\left\|\Delta\mathbf{v}_k\right\| + \left\|\mathbf{v}_{Np}^\infty\right\| \tag{13}$$

or, if the launch velocity is fixed to a given maximum value, as:

$$F(\mathbf{y}) = \sum_{k=1}^{N_p-1}\left\|\Delta\mathbf{v}_k\right\| + \left\|\mathbf{v}_{Np}^\infty\right\| \tag{14}$$

Once the problem in Eq. (12) is formulated for a given sequence of gravity-assist bodies a search over a large range of launch dates, encounter times and gravity-assist characteristics can be started, in order to locate the best gravity-assist trajectory to reach a target planet. Furthermore, if the sequence of planetary encounters is left free, problem in Eq. (12) has to be extended to handle mixed integer-real parameters (MINLP problem). In this case fixing the maximum number of encounters the solution vector becomes:

$$\mathbf{y} = [p_1,...,p_k,...,p_{N_p-1}, v_0^\infty, \alpha, \delta, t_0, T_1, \varepsilon_1, \eta_1, h_1, ...T_k, \varepsilon_k, \eta_k, h_k, ...T_{Np-1}, \varepsilon_{Np-1}]^T \tag{15}$$

where $p_k$ is the identification reference number of planet $k$-$th$ or of any other object available in the database with the orbital elements. In this case the integer part of the problem represents a combinatorial problem in which only few optimal sequences are of interest. It should be noted that if all the sequences of any length are of interest, this is not a clever way of formulating the problem. A better approach would be to decouple the combinatorial part from the real-value part. For the same reason we have not inserted the length of the sequence into the vector **y** since this parameter belongs to a one-dimensional limited and countable set containing just a limited number of interesting elements.

## IV.   Global Search Through a Hybrid Evolutionary Algorithm

### A. Search and Optimization Issues

If $N_P$ is the number of planets in the sequence under investigation, then the number of design parameters increases as $n=4N_p-2$, with a typical number of design parameters ranging from 18 to 22 in the case of a mission with two or three flybys. Supposing the implementation of a systematic search method, one option could be to discretize the search domain and evaluate each possible solution. If only two samples are taken for each variable in the case of a transfer with four flybys, such as the Cassini trajectory, at least $2^{22}$=4,194,304 evaluations of the objective function would be necessary.

On the other hand, if we require to this simple systematic search algorithm to be complete, i.e. to guarantee convergence to a globally optimal solution, for a block-box real-valued problem, an infinite number of samples densely distributed over the whole solution space would be required. Note that completeness of the search and exhaustiveness of the search, in this case, coincide since the algorithm does not have simply to find a solution but has to prove that there are no other better solutions. This makes a complete search practically impossible without any additional information on the landscape of the search space.

If completeness is dropped in practice, then the aim could be to find a good set of solutions in a realistic computational time, while retaining the theoretical completeness of the search algorithm. To this aim we could use a stochastic sampling algorithm based on some efficient heuristics that drives the sampling procedure, such as Evolutionary Algorithms.

Since the probability of finding a solution is dependent on the density of samples in a given region, it would be desirable to prune the unpromising areas of the solution space. This is achieved through a proper blending of an evolutionary based algorithm with a systematic branching strategy, as will be described in the following chapters.

**B. The Evolutionary-Branching Principle**

Evolutionary-Branching (EB) is a hybrid deterministic-stochastic approach to the solution and characterization of constrained and unconstrained multimodal, multivariate nonlinear programming problems with mixed integer-real variables and discontinuous quantities. The EB approach is based on the following principal ideas:

- An evolutionary algorithm is used to explore the solution space $D$, then a branching scheme, dependent on the findings of the evolutionary step, is used to partition the solution domain into subdomains $D_l$. If necessary, on each $D_l$ a new evolutionary search is performed and the process is repeated until a number of good minima and eventually the global one are found (practical convergence) or the entire domain has been partitioned into infinitesimal subdomains (theoretical convergence).

- The search is performed by a population of individuals $\mathbf{y}$, each one represented by a string of length $n$, containing in the first $m$ components integer values and in the remaining $n$-$m$ components real values. A hypercube $\mathbf{S}$ enclosing a region of the solution space surrounding each individual, is then associated to $\mathbf{y}$.

The solution space is explored locally by acquiring information about the landscape within each region **S** and globally by a portion of the population, which is continuously regenerated.

- Each individual can communicate its findings to the others in order to evolve the entire population towards a better status.
- During the evolutionary step a discoveries-resources balance is maintained: a level of resources is associated to each individual and is reduced or increased depending on the number of good findings of the individual.

This particular hybridization of evolutionary algorithms is motivated by the fact that common EA need to be run several times on this problem to provide reliable results. The deterministic step is then intended to reduce the search space at every new run. A further improvement over standard EA has been achieved by an increase of the local search capabilities of each individual. This was obtained by the introduction of a novel mechanism, called perception in the following. A comparison of this approach with known global optimizers is outside the scope of this paper and can be found in Vasile[20] and Di Lizia et al.[21].

### C. Environment Perception

Each region **S** is evaluated using a mechanism called *perception*. This operator samples the environment in order to improve the status of the individuals. A new region **S** is then associated to the best discovered value, resulting in a migration of the individual towards a place where better resources are expected. For this reason each hypercube **S** is here called *migration region*. The set of samples in **S** is generated with the following procedure: a first sample $\mathbf{y}^{(1)}$ is generated, by mutation of **y** (mutation operators are taken from standard real-valued Evolutionary Algorithms sampling mechanisms[24]), then a linear extrapolation is performed. Extrapolation generates a new sample $\mathbf{y}^{(2)}$ on the side of the best one between $\mathbf{y}^{(1)}$ and **y** as follows:

$$\mathbf{y}^{(2)} = \nu(\mathbf{y}^{(1)} - \mathbf{y}) + \mathbf{y}^{(1)} \tag{16}$$

where $\mathbf{y}^{(1)}$ is here assumed to be better than **y** and $\nu$ is a random number taken from a uniform distribution. The two resulting individuals $\mathbf{y}^{(1)}$, $\mathbf{y}^{(2)}$ and the parent **y** are then used to generate a third individual by using second order interpolation mating. If **p** is the vector difference between **y** and $\mathbf{y}^{(2)}$ and $f, f^{(1)}, f^{(2)}$ are the fitness values for the three individuals $\mathbf{y}, \mathbf{y}^{(1)}, \mathbf{y}^{(2)}$ respectively, then second order interpolation mating generates an individual building a one-dimensional model of the fitness function and taking as a new individual the minimum of the resulting:

$$\mathbf{y}^{(3)} = \mathbf{y} + \mathbf{p}\chi_{min}, \tag{17}$$

$$f_{min} = a(\mathbf{y}, \mathbf{y}^{(1)}, \mathbf{y}^{(2)})\chi_{min}^2 + b(\mathbf{y}, \mathbf{y}^{(1)}, \mathbf{y}^{(2)})\chi_{min} + f(\mathbf{y}) \tag{18}$$

The procedure is repeated until either an improvement of **y** is found or a number of samples equal to the level of resources has been generated. The level of resources is increased by one unit if an improvement occurs and is decreased by one unit if nothing is found. The upper limit has been fixed to the number of coordinates and the lower limit set to 1.

The contraction or expansion of each region **S** is regulated through a parameter $\rho$, called migration radius, which depends on the findings of the perception mechanism. The migration radius is defined as the ratio between the value of the distance from the boundary $\mathbf{b}_j$ of the migration region of the j-th individual and the value of the distance from the corresponding boundary **b** of the domain $D$:

$$\rho_j = \frac{b_{i,j} - y_{i,j}}{b_i - y_{i,j}} \tag{19}$$

If none of the samples in **S** is better than **y** the radius is reduced according to:

$$\rho_j = \begin{cases} \max([1e-8, \delta y_{min}]) & \text{if } \delta y_{min} \geq \varepsilon \rho_j \\ \varepsilon \rho_j & \text{if } \delta y_{min} < \varepsilon \rho_j \end{cases} \tag{20}$$

where $\varepsilon$ has been set to 0.5 and $\delta y_{min}$ is the distance of the best sample $\mathbf{y}^*$, among the ones in the migration region, from the $\mathbf{y}_j$, normalized with respect to the dimensions of the migration region:

$$\delta y_{min} = \sqrt{\sum_{i=1}^{n} \left(\frac{y_i^* - y_{i,j}}{S_{i,j}}\right)^2} \tag{21}$$

where, for j-th individual and for dimension $i$, $S_{i,j}$ is the difference between the value of the upper bound and of the lower bound and the summation is over non-zero dimensions. Now, if from generation $k$ to generation $k+1$ the differential improvement $\Delta f_j$ (i.e. the difference between the function $f_j$ at generation $k$ minus $f_j$ at generation $k+1$) increases, the migration radius is recomputed according to the prediction:

$$\rho_j^{(k+1)} = \rho_j^{(k)} \theta \log(e - 1 + j) \tag{22}$$

where $\theta$ is equal to 2 in this implementation.

**D. Communication Mechanisms**

At the end of a full evolution step those individuals that have improved their status are inserted in a communication list and exchange information with an equal number of randomly selected individuals from the entire population. The individuals can communicate through a simple exchange of their components, the linear extrapolation in Eq. (16) or through a linear interpolation operator given by:

$$\mathbf{y}^{(2)} = \nu(\mathbf{y}^{(1)} - \mathbf{y}) + \mathbf{y} \tag{23}$$

The interpolation operator is used also to prevent crowding of more than one individual in the basin of attraction of the same solution: if the reciprocal distance among two or more individuals falls down below a given threshold, the worst one is mated with the boundaries of the subdomain $D_l$ thus projecting the individual into a random point within $D_l$, according to the following relation:

$$y_i^{(1)} = \nu b_i + (1-\nu) y_i \tag{24}$$

**E. Ranking**

At each evolutionary step the entire population of *npop* individuals is ranked from the best to the worst and the best $n_e$ individuals are allowed to use the perception mechanism while the others are either hibernated (i.e. no operator is applied) or mutated. The probability of being mutated or hibernated depends on their ranking: the lower the rank position is, the higher the chance to be mutated.

**F. Branching Step**

Even though the evolutionary step can find several optima and eventually the global one, convergence is not guaranteed due to the stochastic nature of the process. Therefore, a systematic decomposition of the solution space is performed on the basis of the output of the evolutionary algorithm.

The initial domain *D* is recursively partitioned generating a number of smaller and smaller subdomains $D_l$. The partitioning, or branching, process selects the worst individual $\mathbf{y}_{worst}$, found by the evolutionary step within $D_l$ (with $D_0 = D$), and decomposes $D_l$ into *2n* subdomains, or nodes, by cutting the coordinates along the components of $\mathbf{y}_{worst}$ (a safeguard mechanism prevents cuts too close to a boundary by moving the cutting point to the middle of the interval).

For each node the ratio between the relative number of individuals and the relative volume is computed and the resulting quantity defines how necessary a further exploration of the node is:

$$\varpi_{D_q} = \frac{\sum_{D_q} j}{\sum_{D_l} j} \Bigg/ \sqrt[n]{\frac{V_{D_q}}{V_{D_l}}} \qquad q=1,\ldots,2n \qquad (25)$$

where the volumes $V_{D_q}$ and $V_{D_l}$ are computed taking only edges with a non-zero dimension. This quantity is then added to a fitness $\varphi_{D_q}$ defined as:

$$\varphi_{D_q} = \begin{cases} \dfrac{\dfrac{1}{J}\sum_{j=1}^{J} f_j - f_{best}}{f_{worst} - f_{best}} & \text{if } J \neq 0 \\ 1 & \text{otherwise} \end{cases} \qquad (26)$$

where $J$ is the number of individuals in domain $D_q$ while $f_{best}$ and $f_{worst}$ are respectively the best fitness values in the whole population and the worst fitness value in the whole population. The node is then qualified by the quantity:

$$\psi_{D_q} = \sigma \varpi_{D_q} + (1-\sigma)\varphi_{D_q} \qquad (27)$$

where $\sigma$ is the weighting factor that weights how reliable the result coming from the evolution step is considered. If $\sigma$ is 0, only the nodes with low fitness are explored because the evolutionary algorithm is considered reliable enough to explore exhaustively the domain $D_l$ without leaving any region unexplored. On the other hand if $\sigma$ is 1 the result from the evolutionary algorithm is considered to be not reliable due to a premature convergence or to a poor exploration of the solution space.

Among all $2n$ nodes only the best pair according to quantity in Eq. (27) are selected, then, in order to avoid the rediscovery of already found minima, a second cut is performed at the best converged individual $\mathbf{y}_{best}$ and two additional nodes are selected. Therefore at end of the branching step the subdomain $D_l$ has been divided into three new subdomains. These nodes are added to the list of all the $L$ potentially interesting subdomain such that:

$$D = \bigcup_{l=1}^{L} D_l \qquad (28)$$

Once the list is ranked from the best to the worst, according to Eq. (27), the best node is selected for further exploration.

.

## G. Stopping Criteria

In this work, three stopping criteria are used: the maximum number of function evaluations, the number of times subdomains have been branched without improvement, or branching level, and the convergence of the best individual in the filter (only for the evolutionary step). All of them are based on some heuristics and not on any rigorous proof of global convergence. It should be noted, however, that the branching scheme is devised to asymphtotically partition the whole domain $D$ into infinitesimal subdomains and therefore to converge globally. Finally, solution accuracy is enhanced starting a local search with an SQP algorithm from each one of the best solutions found by the evolutionary-branching algorithm. It should however be pointed out that, in the cases of a direct transfer or single flyby trajectory, the value of the solution obtained by evolutionary-branching scheme, within a typical number of function evaluations, could not be improved further by the SQP algorithm since local convergence had already been achieved. Pseudo-codes both for the branching algorithm and for the evolutionary algorithm are represented in Fig. 7.

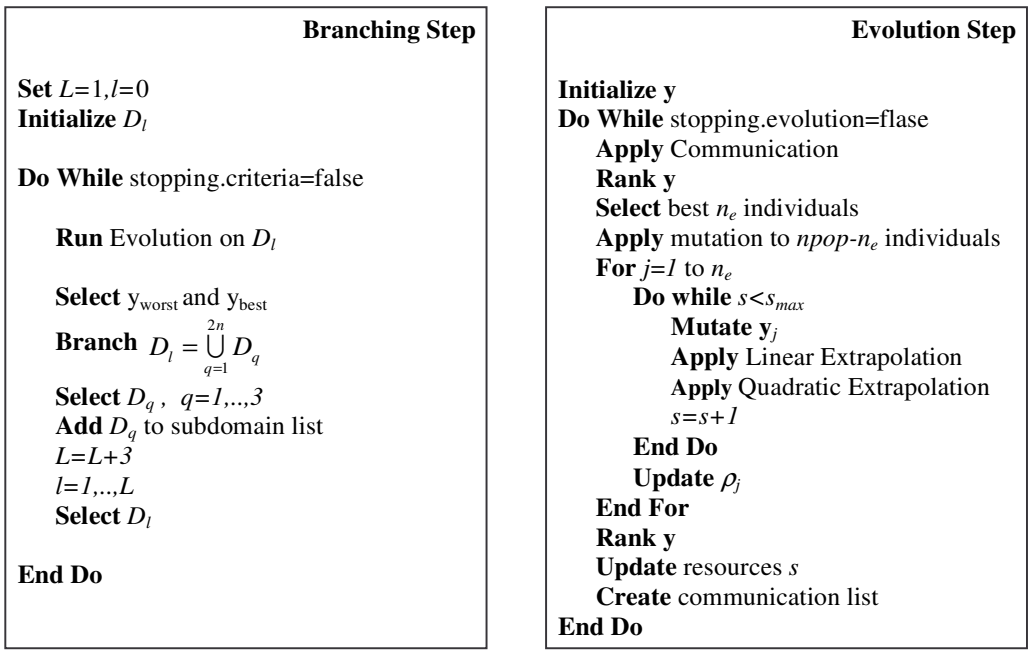

**Fig. 7 Pseudo-code for the evolutionary-branching algorithm.**

## V. Case Studies

The methodology presented in the previous chapters is at the core of a preliminary design tool, implemented in Matlab, called IMAGO (Interplanetary Mission Analysis Global Optimization). IMAGO combines the described trajectory model with the evolutionary-branching algorithm called EPIC. In the following the tool is tested on a number of difficult cases: the design of multiple transfer options to Jupiter in the case of free flybys sequence, the design of the Cassini trajectory, the design of a NEO interception mission, and finally the design of the Rosetta mission. All the tests have been performed on a 2GHz Intel Centrino processor with Windows XP.

**A. Missions to Jupiter**

This first case study requires the design of optimal free-sequence transfers to Jupiter: only the arrival planet, the departure planet and the maximum number of swingbys have been defined, while the sequence of planetary encounters is left free. Table 1 shows the solution domain for the Earth-Jupiter transfer problem, where the pericenter altitudes $\bar{h}_k$ are normalized with respect to the planet mean radius.

**Table 1 Solution domain for the Earth Jupiter problem.**

| Bounds | $p_1$ | $p_2$ | $p_3$ | $t_0$ (MJD2000) | $T_1$ (day) | $\bar{h}_1$ | $T_2$ (day) | $\varepsilon_1$ | $\bar{h}_2$ | $T_3$ (day) | $\varepsilon_2$ | $\bar{h}_3$ | $T_4$ (day) | $\varepsilon_3$ |
|---|---|---|---|---|---|---|---|---|---|---|---|---|---|---|
| $U_b$ | 2 | 2 | 2 | 3650 | 80 | 0.04 | 80 | 0.01 | 0.04 | 80 | 0.01 | 0.04 | 600 | 0.01 |
| $L_b$ | 4 | 4 | 5 | 7300 | 400 | 14 | 800 | 0.9 | 14 | 800 | 0.90 | 14 | 2000 | 0.90 |

The number of possible encounters has been limited to three, since a larger number would excessively increase the total transfer time while adding little benefits. The launch date can vary between 2010 and 2020 in order to look for multiple alternative launch opportunities. The range of possible planets goes from Venus (number 2) to Mars (number 4) for the first two gravity maneuvers, and from Venus to Jupiter (number 5) for the third one. It can be seen that the upper bound on the last swing-by, allows flying directly to Jupiter and thus to investigate solutions with a reduced number of gravity maneuvers. For this first test case a population of 40 individuals has been used with a filter dimension of 20 individuals and 4 levels of branching. This resulted in eight hours of computing time and about 1000 solutions with different combinations of planets and launch dates. For a more extensive search the algorithm was run three times and we reported in Table 2 some interesting solutions that represent a valid alternative to the trajectories reported in the work of Petropopulos et al.[13].

**Table 2 Summary of some interesting options for the transfer to Jupiter.**

| Sequence | Launch date (DD/MM/YY) | $v_0^\infty$ (km/s) | $\Delta v$ (km/s)[a] | TOF (day) | Arrival $v^\infty$ (km/s) |
|---|---|---|---|---|---|
| EVEEJ | 17/08/10 | 2.90 | 2.103 | 1896 | 5.640 |
| EVEEJ | 01/04/12 | 3.26 | 0.0 | 2706 | 6.450 |
| EVEEJ | 27/10/13 | 3.69 | 0.134 | 2082 | 5.798 |
| EVEEJ | 26/10/13 | 3.68 | 0.174 | 2005 | 6.270 |
| EVVEJ | 28/10/13 | 3.87 | 0.620 | 2404 | 6.305 |
| EMMJ | 20/03/14 | 3.31 | 3.100 | 2243 | 4.640 |
| EVEMJ | 20/01/17 | 3.68 | 0.648(2)+0.544(3) | 3112 | 6.639 |
| EVEEJ | 14/03/20 | 3.09 | 0.109 | 2519 | 6.210 |
| EVEEJ | 20/03/20 | 3.01 | 0.546 | 2740 | 5.564 |
| EVEEJ | 23/03/20 | 2.99 | 0.125 | 2545 | 5.575 |
| EVEEJ | 26/10/21 | 3.34 | 0.223(2)+1.327(3) | 1979 | 5.609 |
| EVVEJ | 11/12/21 | 4.22 | 0.183 | 2577 | 6.102 |
| EVEEJ | 01/06/23 | 3.30 | 0.0 | 2498 | 5.746 |
| EVEEJ | 02/07/26 | 3.29 | 0.853 | 2133 | 6.359 |
| EVEJ | 18/10/29 | 3.13 | 1.668 | 1473 | 6.080 |

[a]The parenthetical numbers indicate the flyby number after which a maneuver occurs

It can be noted that the tool was able to identify the typical optimal sequences for a mission to Jupiter, such as EVVEJ or EVEEJ and to provide for these sequences some interesting solutions in terms of escape and arrival velocity. Furthermore it is interesting to underline the Mars resonant solution; although no particular prescription on the kind of solution was set, and additionally no particular model is implemented in order to generate solutions that exploit resonant swing-bys, the EMMJ transfer, plotted in Fig. 8, was found as a result of the search step. It is an interesting solution since it presents an unusual possibility for a transfer to Jupiter that, although requiring a large deep space correction, could be effectively implemented through a low-thrust propulsion system.

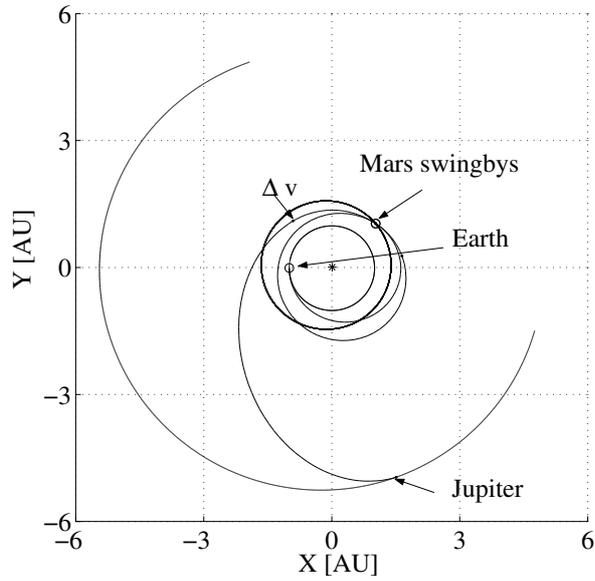

**Fig. 8 Projection into the ecliptic plane of the sequence EMMJ.**

**Table 3 Summary and comparison of different options for the 2009 launch window.**

| Solutions | Launch date (DD/MM/YY) | $v_0^\infty$ (km/s) | $\Delta v$ (km/s) | TOF (day) | Arrival $v^\infty$ (km/s) |
|---|---|---|---|---|---|
| 2009 ESOC | 04-24/03/09 | 4.08-4.38 | 0.05 | 2046-2096 | 5.87-5.99 |
| Petropoulos et al. [13] | 21/01/09 | 3.68 | 0.00 | 2299 | 5.50 |
| Sol.1 IMAGO | 15/03/09 | 4.04 | 0.00 | 2067 | 6.00 |
| Sol.2 IMAGO | 17/03/09 | 4.09 | 0.00 | 2209 | 5.55 |
| Sol.3 IMAGO | 31/01/09 | 3.69 | 0.03 | 2238 | 5.59 |

In order to extend the validation of the tool, the sequence EVEEJ has been further investigated in the year 2009 and 2010 and compared to some solutions found in the literature[13,25]. Transfer time to Venus and to the first Earth flyby can vary between 70 and 550 days, while the transfer time to the second earth flyby can be as long as 800 days, thus allowing the possibility of having a two years resonant orbit after the first flyby. Lastly the transfer time to Jupiter can vary between 400 and 1600 days. A deep space maneuver can be applied up to 90% of the transfer time after each flyby, while the maximum altitude of the swingby maneuver is bounded below by 300 km for Venus and by 1000 km for Earth.

**Table 4 Summary and comparison of different options for the 2010 launch window**

| Solutions | Launch date (DD/MM/YY) | $v_0^\infty$ (km/s) | $\Delta v$ (km/s) | TOF (day) | Arrival $v^\infty$ (km/s) |
|---|---|---|---|---|---|
| 2010 ESOC | 01/07/10-31/08/10 | 2.81-2.97 | 0.500 | 2379-2426 | 5.92-6.33 |
| Petropulos et al.[13] | 21/07/10 | 2.94 | 0.360 | 2409 | 6.21 |
| Sol.1 IMAGO | 04/08/10 | 2.73 | 0.380 | 2387 | 5.53 |
| Sol.2 IMAGO | 03/08/10 | 2.91 | 0.538 | 2424 | 5.80 |

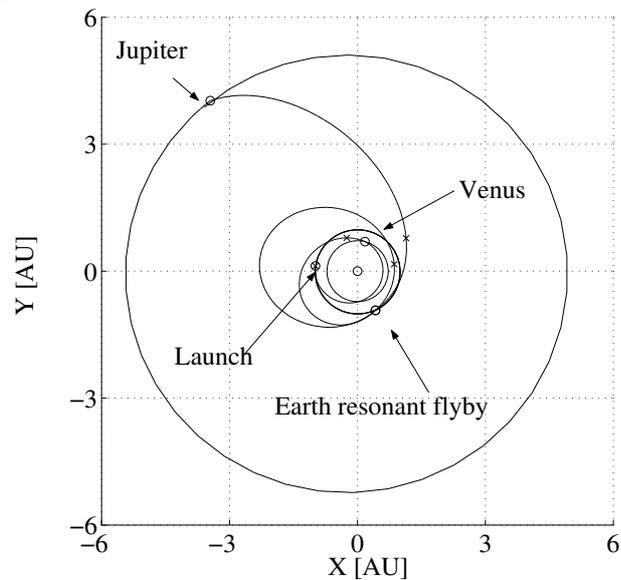

**Fig. 9 Solution for the EVEEJ sequence for 2009 launch.**

For this test, which has a fixed sequence and smaller bounds on the launch date, it was sufficient to use only the evolutionary step. A population of 40 individuals has been evolved for a maximum number of 200,000 function evaluations. The filter dimension has been set equal to 25, thus favoring convergence rather than exploration. Although not exploiting the whole investigation capability offered by the search algorithm, as can be seen in Table 3, there is a good agreement between the solution computed by IMAGO and the reference one. In Fig. 9 and Fig.10 two examples for the sequence EVEEJ for the 2009 and 2010 launch date are shown, while Table 3 and 4 show respectively three and two solutions found for the 2009 and 2010 launch date and a comparison with the best solutions found by Petropoulos et al.[13] for the same launch year.

It should be noted that for both launch options, IMAGO was able to locate, in approximately twenty minutes of computational time, within the launch window some solutions already known in literature. However they show an improvement over the reference solution both in terms of total Δv and transfer time.

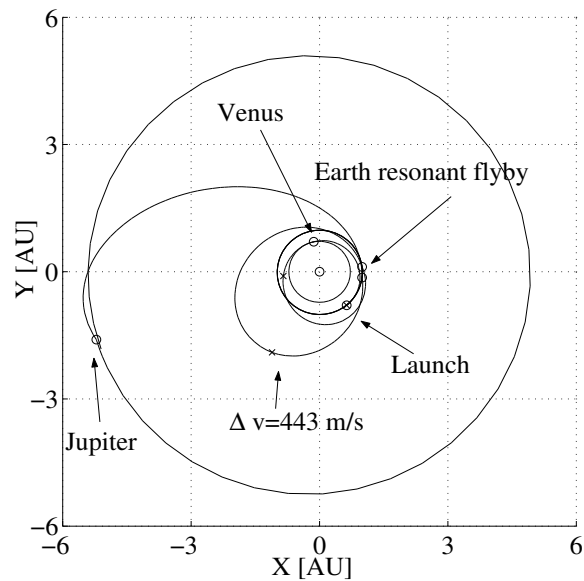

**Fig. 10 Solution for the EVEEJ sequence for 2010 launch.**

**B. Cassini Trajectory Design**

The design of the Cassini trajectory is an interesting test case since it is one of the most complex MGA trajectories designed for a mission to an outer planet. Furthermore the arrival velocity must be kept sufficiently low in order to allow the insertion into the Saturn planetary system and the high mass budget of the spacecraft limits the launch hyperbolic escape velocity provided by the launcher, below 4 km/s. The launch date investigated is 1997 and the other bounds on the design parameters are reported in Table 5; a deep space maneuver can occur up to 90% of the transfer time to each planet.

As can be seen in Table 5 the bounds on the transfer time are quite broad, since in every preliminary design phase there is a lack of preliminary information on possible good options, and every possibility should be investigated. This test case thus reproduces the typical preliminary design conditions when little information is available. The preliminary design is then initiated by looking for solutions minimizing the function in Eq. (13). For this test a population of 40 individuals has been used with a filter dimension equal to 25 and two levels of branching. The total number of function evaluation was a little bit more than 800,000 with a final number of around 200 generated solutions. Fig.11 shows all the local minima that have been found in the search domain: most of them are distributed around the date -800 MJD which is the optimal launch window which collects all the solutions with a total Δv of approximately 10 km/s (the total Δv of Cassini). Fig. 12 shows the EVVE part of the transfer with the Venus-targeting Δv maneuver.

**Table 5 Bounds on the design parameters for the Cassini trajectory.**

| Bounds | $T_1$ (day) | $T_2$ (day) | $T_3$ (day) | $T_4$ (day) | $T_5$ (day) | $\bar{h}_1$ | $\bar{h}_2$ | $\bar{h}_3$ | $\bar{h}_4$ |
|---|---|---|---|---|---|---|---|---|---|
| $L_b$ | 100 | 100 | 30 | 400 | 800 | 0.05 | 0.05 | 0.15 | 0.7 |
| $U_b$ | 400 | 500 | 300 | 1600 | 2200 | 5 | 5 | 5.5 | 290 |

In Table 6 the reference solution for the Cassini mission and two of the best solutions found by IMAGO have been compared. It is interesting to highlight the similarity of the first solution with the reference solution both in terms of launch dates and total Δv, while the second solution presents an interesting saving in the arrival velocity at the expense of a higher transfer time.

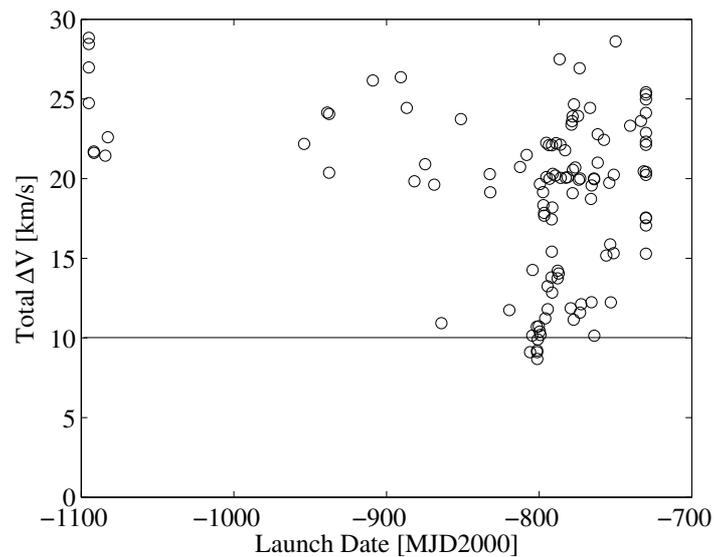

**Fig. 11 Distribution of the optimal solutions found for a 1997 launch for Cassini.**

**Table 6 Summary and comparison of different options for Cassini.**

| Mission parameter | | CASSINI | Sol 1 IMAGO | Sol 2 IMAGO |
|---|---|---|---|---|
| Launch date (DD/MM/YY) | | 15/10/97 | 20/10/97 | 17/10/97 |
| $v_0^\infty$ | (km/s) | 3.93 | 4.04 | 4.03 |
| E-V TOF | (day) | 194 | 191 | 191 |
| V-V DSM | (km/s) | 0.471 | 0.432 | 0.414 |
| V-V TOF | (day) | 425 | 421 | 420 |
| V-E DSM | (km/s) | 0 | 0 | 0 |
| V-E TOF | (day) | 54 | 53 | 53 |
| E-J DSM | (km/s) | 0 | 0.132 | 0 |
| E-J TOF | (day) | 499 | 493 | 540 |
| J-S DSM | (km/s) | 0.376 | 0 | 0 |
| J-S TOF | (day) | 1267 | 1216 | 1656 |
| Total $\Delta v$ | (km/s) | 10.14 | 10.18 | 9.06 |

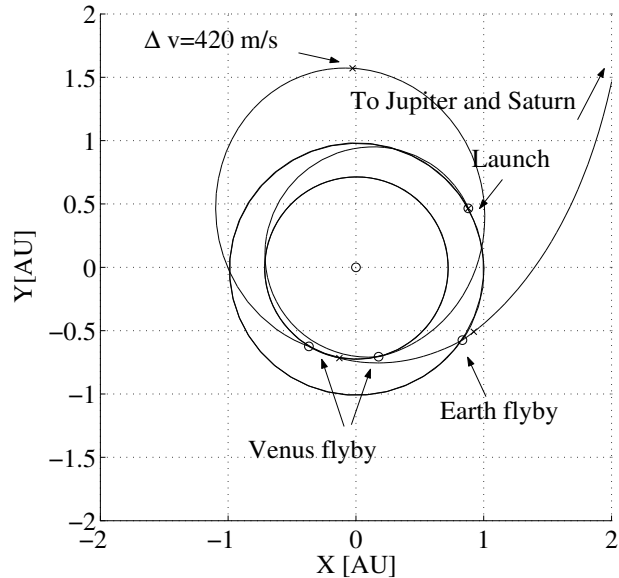

**Fig. 12 Close view of the sequence of the flybys of the inner planets for Cassini.**

As further verification the same problem with the same model has been solved with a simple multistart technique: 100 samples of each parameter have been taken with a Latin Hypercube distribution and the best three samples have been optimized with a local procedure using the Matlab function fmincon. The procedure has been repeated for 30 independent runs for a total of about 1,600,000 function evaluations. Though the overall number of function evaluations is twice that of IMAGO, the total number of good solutions is lower. Furthermore, though this procedure was run several times, the best solutions obtained with the multistart approach were worse than the ones produced by IMAGO. This proves that a simple random approach with a simple heuristics cannot solve efficiently this problem. This was proven to be even more true for problems of higher complexity or for an extended launch window.

**C. Missions to Asteroid (10302) -1989ML**

Missions to asteroids have been gaining increasing interest over the last few years, both for scientific reasons and for hazard mitigation. In particular the mission Don Quijote,[26] currently under study in Europe, is planning to fly two twin spacecraft to asteroid (10302)-1989ML in order to perform scientific investigations with one satellite, while the other impacts the asteroid at a high speed (greater than 10 km/s). The trajectory of the impacting spacecraft has been taken as the reference solution for this analysis.

The search was performed over a launch window in 2011 and only the solutions with an arrival velocity greater than 10 km/s have been retained. In this case, only the evolutionary step was used with a population of 30 individuals and a maximum number of 100,000 function evaluations.

The reference solution proposes a first Earth swing-by after 180 days in order to separate the two spacecraft that have been launched together, therefore the selected sequence is EEV-asteroid. Here the bounds on the E-E transfer orbit (Table 7 summarizes the bounds on all the design parameters) can range from 80 to 400 days hence allowing the investigation of other possible options and particularly a one-to-one resonant flyby with the Earth (see Table 8). Transfers with a 180 days Earth-to-Earth flyby are the most recurrent solutions and can be classified into two groups: short period Earth-Venus, Venus-Asteroid transfers (see Fig. 13) and long period Earth-Venus, Venus-Asteroid transfers (see Fig. 14). The latter have a bound orbit before the encounter with Venus and a second bound orbit before the impact.

**Table 7 Bounds on the design parameters for the mission to asteroid (10302) -1989ML.**

| Bounds | $t_0$(MJD2000) | $\varepsilon_1$ | $\bar{h}_1$ | $T_1$(days) | $\varepsilon_2$ | $\bar{h}_2$ | $T_2$(days) | $\varepsilon_3$ | $T_3$(days) |
|--------|----------------|-----------------|-------------|-------------|-----------------|-------------|-------------|-----------------|-------------|
| $U_b$ | 4015 | 0.01 | 0.1 | 80 | 0.01 | 0.1 | 80 | 0.01 | 600 |
| $L_b$ | 5100 | 0.90 | 14 | 400 | 0.90 | 14 | 800 | 0.90 | 2000 |

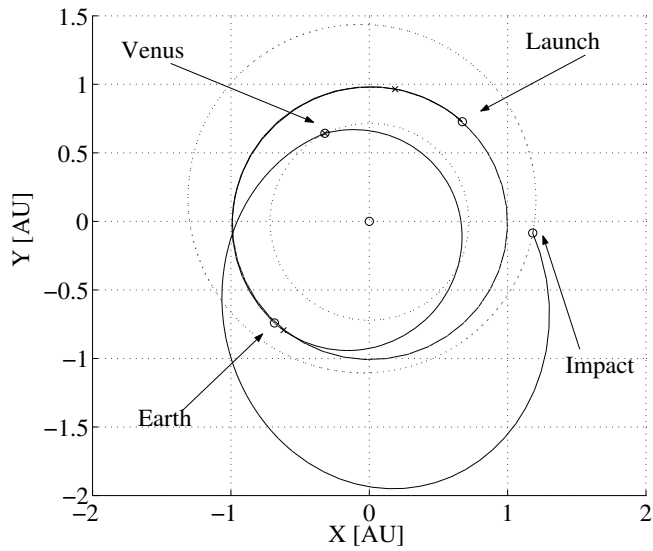

**Fig. 13 Projection into the ecliptic plane of Sol 3.**

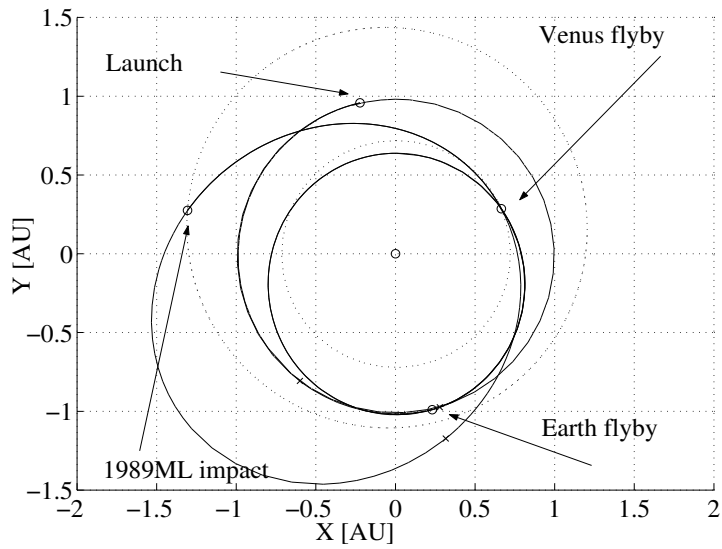

**Fig. 14 Projection into the ecliptic plane of Sol 1.**

**Table 8 Summary of different solutions for the mission to (10302) -1989ML.**

| Mission parameter | Sol 1 IMAGO | Sol 2 IMAGO | Sol 3 IMAGO | Sol 4 IMAGO | Sol 5 IMAGO |
|---|---|---|---|---|---|
| Launch date (DD/MM/YY) | 05/01/11 | 01/01/11 | 03/11/11 | 19/10/11 | 26/03/11 |
| $v_0^\infty$ (km/s) | 3.99 | 3.59 | 3.85 | 3.64 | 3.8 |
| E-E DSM (km/s) | 0 | 0 | 0 | 0 | 0 |
| E-E TOF (day) | 182 | 182 | 179 | 178 | 365 |
| E-V DSM (km/s) | 0 | 0.314 | 0 | 0.341 | 0.327 |
| E-V TOF (day) | 641 | 726 | 168 | 179 | 191 |
| E-A DSM (km/s) | 0 | 0 | 0.389 | 0 | 0 |
| E-A TOF (day) | 603 | 791 | 478 | 471 | 469 |
| Arrival $v^\infty$ (km/s) | 13.08 | 13.15 | 16.81 | 16.34 | 15.77 |

**D. Rosetta Trajectory Design**

The preliminary design of the Rosetta mission represents an interesting test case for two different reasons. Firstly the trajectory rendezvous with a comet on an elliptical orbit, after a number of swing-bys and deep space maneuvers. Secondly it is an interesting example of the demand for a quick investigation of different back up options. The Rosetta mission, originally designed to rendezvous with comet Wirtanen, was not launched in January 2003 due to technical problems with Ariane-5 and was completely re-designed for a new target and set to launch in 2004. This required a new investigation and characterization of all the trajectory opportunities, alongside the definition of a new target. For this reason this test case has been approached supposing very little knowledge on the structure of the solutions in terms of launch date and deep space corrections. The sequence has been chosen to be EEMEE-comet which is the solution adopted for the new Rosetta mission.[27]

The design has been conducted in two different steps; as a first step a deep investigation, for a launch date in 2004 with a fixed velocity of 3.546 km/s due to the launcher capability, has been performed using the large bounds on the design parameters reported in Table 9. This is the typical situation that mission analysts would face during the preliminary design of a new mission, in order to assess good sequences. For this first investigation the branching procedure was largely exploited, using two levels of branching and a maximum number of 200,000 evaluations for each evolutionary step, with a population of 40 individuals. This parameter setting results in a maximum number of some 800,000 total evaluations, but yields an extensive characterization of the possible solutions and the local optima in the domains.

**Table 9 Bounds on the design parameters for the Rosetta trajectory.**

| Bounds | $t_0$ (MJD2000) | $T_1$ (day) | $\varepsilon_1$ | $\bar{h}_1$ | $T_2$ (day) | $\varepsilon_2$ | $\bar{h}_2$ | $T_3$ (day) | $\varepsilon_3$ | $\bar{h}_3$ | $T_4$ (day) | $\varepsilon_4$ | $\bar{h}_4$ | $\varepsilon_5$ | $T_5$ (day) |
|---|---|---|---|---|---|---|---|---|---|---|---|---|---|---|---|
| $U_b$ | 1460 | 300 | 0 | 0.06 | 150 | 0 | 0.06 | 150 | 0 | 0.06 | 300 | 0 | 0.06 | 0 | 700 |
| $L_b$ | 1825 | 400 | 0.9 | 9 | 800 | 0.9 | 9 | 800 | 0.9 | 9 | 800 | 0.9 | 9 | 0.9 | 1850 |

**Table 10 Summary of different options for the Rosetta mission.**

| Mission parameters | | Rosetta | IMAGO: Sol 1 | IMAGO: Sol 2 | IMAGO: Sol 3 | IMAGO: Sol 4 |
|---|---|---|---|---|---|---|
| Launch date (DD/MM/YY) | | 02/03/04 | 05/03/04 | 07/03/04 | 27/02/04 | 12/03/04 |
| $v_0^\infty$ | (km/s) | 3.546 | 3.546 | 3.546 | 3.546 | 3.546 |
| E-E TOF | (day) | 367 | 361 | 363 | 369 | 370 |
| E-M TOF | (day) | 722 | 722 | 719 | 721 | 728 |
| M-E TOF | (day) | 261 | 263 | 263 | 262 | 256 |
| E-E TOF | (day) | 730 | 729 | 730 | 730 | 730 |
| E-P67 TOF | (day) | 1650 | 1825 | 1825 | 1642 | 1640 |
| Total $\Delta v$ | (km/s) | 1.700 | 1.763 | 1.581 | 1.702 | 1.921 |

The results of this first extensive search are reported in Fig.15 showing all the local optimal solutions located by the algorithm in a computational time of 2.5 hours. It can be noticed that there are two groups of good solutions for a launch date around 1455 MJD and a launch date around 1730 MJD which present a total $\Delta v$ of some 2-2.5 km/s and there are some other interesting options although less crowded for a launch around the 1530 MJD. Starting from this information the bounds on the search have been reduced around these groups of solutions and a further investigation has been performed producing a large variety of good different solutions comparable with the reference one adopted for the actual mission. In this case only the evolutionary step has been used with a maximum number of 100,000 evaluations, yielding some very good solutions, some of which are presented in Table 10, where the characteristics of the solution found are reported, while Fig. 16 shows the most similar to the reference one.

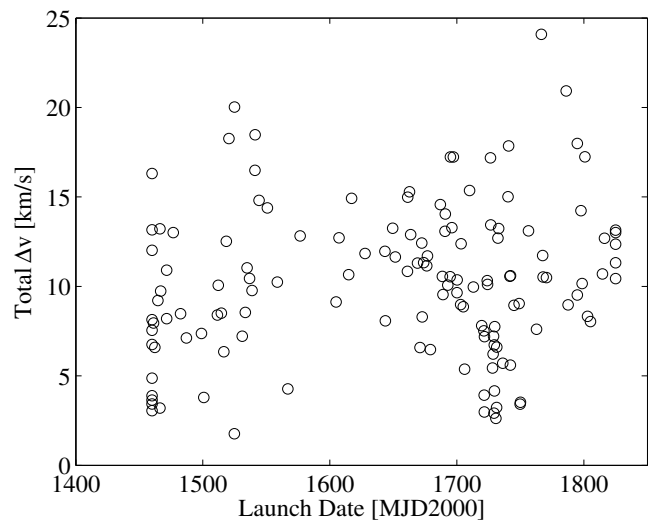

**Fig. 15 Distribution of the optimal solutions for the Rosetta preliminary design.**

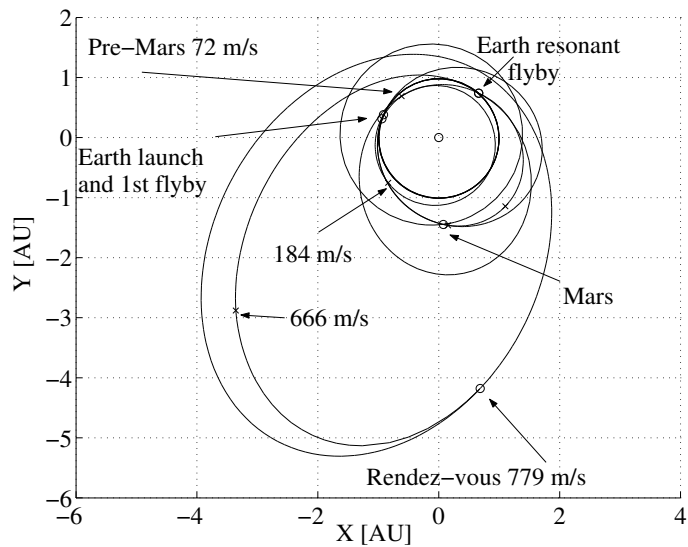

**Fig. 16 Projection into the ecliptic plane of Sol 3.**

## VI. Conclusions

This paper has presented a novel approach to the automatic preliminary investigation of multiple gravity-assist interplanetary trajectories. An initial analysis on the relation between model completeness and solution search complexity has led to the development of a simple trajectory model sufficiently accurate for preliminary design and sufficiently complete to preserve important classes of possible solutions. The analysis performed on the domain structure of each simple sub-problem, composing the whole model, although not exhaustive, has revealed a structure of the search spaces that is generally multimodal and non-convex. These characteristics are a direct consequence of the model complexity and could be exploited to improve the search for a solution. In this paper a general global optimization approach, not problem dependent, is proposed. This novel approach blends together an evolutionary based algorithm with a deterministic branching strategy that is used to partition the whole solution space in subdomains.

The proposed methodology has been tested on a suite of typical mission design problems of high complexity and has proved to yield good results both in terms of quality of the solutions and quantity of information they provide. In particular the nominal Cassini trajectory has been correctly reconstructed starting from very large bounds on the design parameters and in addition some cheaper options in terms of total $\Delta v$ have been found. The Rosetta and the asteroid impact test cases have shown how the proposed methodology is able to characterize highly complex trajectory problems yielding a good deal of different trajectory options with a total $\Delta v$ comparable to the reference solutions. Finally on the free-sequence problem the search algorithm has performed pretty well providing some interesting families of solutions including a few unexpected ones. Furthermore the comparison with the Cassini and Rosetta reference solutions, which represent the optimal trajectory actually flown, shows how the results obtained with the preliminary analysis of IMAGO, are highly accurate with respect to all the main features of the trajectory.

These results demonstrate that the proposed methodology is able to effectively tackle the complexity of multigravity-assist interplanetary transfers while avoiding being trapped in local minima, even in highly multimodal domains, where a classical optimization approach could fail. Moreover, the quality of the solutions and the low computational time make the approach appealing for the development of preliminary design tools.

Finally, the analysis performed in this paper suggests that the search for a solution can be further improved by reducing the solution space. This can be achieved by solving each sub-problem in an incremental fashion and pruning part of the search domain from the early stages of development of the trajectory even for long sequences of gravity maneuvers with up to 15 or more celestial bodies as in the work of Heaton et al.[15] This is the subject of current research and will be presented in future papers.